\documentclass[a4paper,12pt]{article} 
\usepackage{graphicx,latexsym,amsmath,amsthm,amssymb,color,indentfirst,xfrac}
 \usepackage[usenames,dvipsnames]{pstricks}
 \usepackage{epsfig}
\usepackage{pst-grad} 
 \usepackage{pst-plot} 

\newtheorem{satz}{thm1}

\newtheorem{sat}{thm2}

\newtheorem{thm}[sat]{Theorem}
\newtheorem{lem}[satz]{Lemma}

\newtheorem{prop}[satz]{Proposition}
\newtheorem{ex}[satz]{Example}
\theoremstyle{definition}

\newcommand{\be}{\begin{equation}}           
\newcommand{\ee}{\end{equation}}

\newcommand{\ba}{\begin{align}}                
\newcommand{\ea}{\end{align}}

\newcommand{\bal}{\begin{align*}}              
\newcommand{\eal}{\end{align*}}

\newcommand{\bxx}{\begin{ex}}

\newcommand{\exx}{\end{ex}}

\newcommand{\txsit}

\newenvironment{pr}
{\begin{trivlist}
\item[\hskip\labelsep{\bf Proof.}]}                     
{$\hfill\Box$\end{trivlist}}

\title{Defective 2-colorings of planar graphs without 4-cycles and 5-cycles}
\author {Pongpat Sittitrai\\ 
{\small\em Department of Mathematics, Faculty of Science, Khon Kaen University, 40002, Thailand }\\  
{\small\em E-mail address: arzakuraln@gmail.com} 
\and Kittikorn Nakprasit \footnote{Corresponding Author} \\ 
{\small\em Department of Mathematics, Faculty of Science, Khon Kaen University, 40002, Thailand }\\
{\small\em E-mail address: kitnak@hotmail.com}}
\date{}

\begin{document}

\maketitle

\begin{center}{\bf Abstract}\end{center}
\indent\indent
Let $G$ be a graph without 4-cycles and 5-cycles. 
We show that the problem to determine whether $G$ is $(0,k)$-colorable 
is NP-complete for each positive integer $k.$  
Moreover, we construct non-$(1,k)$-colorable planar graphs 
without 4-cycles and 5-cycles for each positive integer $k.$
Finally, we prove that $G$ is $(d_1,d_2)$-colorable 
where $(d_1,d_2)=(4,4), (3,5),$ and $(2,9).$

\section{Introduction}
\indent  Let $G$ be a graph with the vertex set $V(G)$ and the edge set $E(G).$  
A \emph{$k$-vertex}, a \emph{$k^+$-vertex}, and  \emph{$k^-$-vertex}  
are a vertex of degree $k$, at least $k$, and at most $k$, respectively. 
The similar notation is applied for faces. 
A \emph{$(d_1,d_2,\dots,d_k)$-face} $f$ is a face of degree $k$ 
where all vertices on $f$ have degree $d_1,d_2,\dots,d_k$. 
If $v$ is not on a 3-face $f$ but $v$ is adjacent to some $3$-vertex on $f,$  
then we call $f$ a \emph{pendant face} of a vertex $v$ and  
$v$ is a \emph{pendant neighbor} of a $3$-vertex $v.$   
A $3$-face (respectively, $2$-vertex)  incident to a $2$-vertex (respectively, $3$-face) 
is called a \emph{bad $3$-face} (respectively, \emph{bad $2$-vertex}). 
Otherwise, it is a \emph{good $3$-face} (respectively, \emph{good $2$-vertex}).

A \emph{$k$-coloring} $c$ (not necessary proper) 
is a function $c:V(G) \rightarrow \{1,\ldots,k\}.$ 
Define $V_i := \{v \in V(G): c(v) = i\}.$  
We call $c$ a \emph{$(d_1,d_2,\dots,d_k)$-coloring} 
if $V_i$ is an empty set or the induced subgraph $G[V_i]$  
has the maximum degree at most $d_i$ for each $i\in \{1,\ldots,k\}$.
A graph $G$ is called  \emph{$(d_1,d_2,\dots,d_k)$-colorable} if 
$G$ admits a $(d_1,d_2,\dots,d_k)$-coloring
Thus the four color theorem \cite{app1},\cite{app2} can be restated as 
every planar graphs is $(0,0,0,0)$-colorable.  
For improper $3$-colorability of planar graph, Cowen, Cowen, and Woodall  
showed that every planar graph is $(2,2,2)$-colorable \cite{c32}.  
Eaton and Hull  \cite{1kk1} proved that $(2,2,2)$-colorability is 
optimal by showing non-$(k,k,1)$-colorable planar graphs for each $k$. 

Gr\"{o}tzsch \cite{g3f} showed that every planar graph without 3-cycles is 
$(0,0,0)$-colorable. The famous Steinberg's conjecture proposes 
that every planar graph without 4-cycles and 5-cycles is also $(0,0,0)$-colorable. 
Recently, this conjecture is disproved by  Cohen-Addad et al \cite{add}. 
One way to relax the conjecture is allowing some color classes to be improper. 
For every planar graph $G$ without  4-cycles and 5-cycles, 
Xu, Miao, and  Wang \cite{xu110} proved that  $G$ is $(1,1,0)$-colorable,  
and Chen et al. \cite{che200}  proved that $G$ is $(2,0,0)$-colorable. 

Many papers investigate $(d_1,d_2)$-coloring of planar graphs in various settings. 
Montassier and  Ochem \cite{mong4} constructed  planar graphs of girth $4$ 
that are not $(i,j)$-colorable for each $i,j$. 
Borodin, Ivanova, Montassier, Ochem, and Raspaud  \cite{bro6} 
constructed planar graphs of girth $6$ that are not $(0,k)$-colorable for each $k.$  
On the other hand, for every planar graph $G$ of girth $5,$ 
Havet and  Seren \cite{h44} showed that  $G$ is  $(2,6)$-colorable 
and $(4,4)$-colorable, and 
Choi and Raspaud \cite{choi35} showed that  $G$ is  $(3,5)$-colorable.  

\indent  Let $G$ be a graph with the vertex set $V(G)$ and the edge set $E(G).$  
A graph $G$ is called  \emph{$(d_1,d_2,\dots,d_k)$-colorable} if $V(G)$ 
can be partitioned into sets $V_1$, $V_2,\dots,$ $V_k$ such that 
the induced subgraph $G[V_i]$ for $i\in[k]$ has the maximum degree at most $d_i$. 
Thus the four color theorem \cite{app1},\cite{app2} can be restated as 
every planar graphs is $(0,0,0,0)$-colorable.  
For improper $3$-colorability of planar graph, Cowen, Cowen, and Woodall  
showed that every planar graph is $(2,2,2)$-colorable \cite{c32}.  
Eaton and Hull \cite{1kk1} and \v{S}krekovski \cite{1kk2} prove that $(2,2,2)$-colorability is 
optimal by showing non-$(k,k,1)$-colorable planar graphs for each $k$. 

Gr\"{o}tzsch \cite{g3f} showed that every planar graph without 3-cycles is 
$(0,0,0)$-colorable. The famous Steinberg's conjecture proposes 
that every planar graph without 4-cycles and 5-cycles is also $(0,0,0)$-colorable. 
Recently, this conjecture is disproved by  Cohen-Addad et al  \cite{add}. 
One way to relax the conjecture is allowing some color classes to be improper. 
For every planar graph $G$ without  4-cycles and 5-cycles, 
Xu, Miao, and  Wang \cite{xu110} proved that  $G$ is $(1,1,0)$-colorable,  
and Chen et al. \cite{che200}  proved that $G$ is $(2,0,0)$-colorable. 

Many papers investigate $(d_1,d_2)$-coloring of planar graphs in various settings. 
For every planar graph $G$ of girth $5,$ 
Havet and  Seren \cite{h44} showed that  $G$ is  $(2,6)$-colorable and $(4,4)$-colorable, and 
Choi and Raspaud \cite{choi35} showed that  $G$ is  $(3,5)$-colorable.  
Borodin, Ivanova, Montassier, Ochem, and Raspaud  \cite{bro6} 
constructed planar graphs of girth $6$ that are not $(0,k)$-colorable for each $k.$  
Montassier and  Ochem \cite{mong4} constructed  planar graphs of girth $4$ 
that are not $(i,j)$-colorable for any $i,j$. 

There are  many papers \cite{bro6,bo2,h44,bo3,bo1,mong4}  
that investigate $(d_1, d_2)$-colorability forgraphs with girth length of $g$  for $g\geq6$; 
see \cite{mong4} for the rich history. 
For example, Borodin, Ivanova, Montassier, Ochem, and Raspaud \cite{bro6} 
constructed a graph in $g_6$ (and thus also in $g_5$) that is not $(0, k)$-colorable for any $k$. 
The question of determining if there exists a finite $k$ where all graphs in $g_5$ 
are $(1, k)$-colorable is not yet known and was explicitly asked in \cite{mong4}. 
On the other hand, Borodin and Kostochka \cite{bo2} and Havet and Sereni \cite{h44}, respectively, 
proved results that imply graphs in $g_5$ are $(2, 6)$-colorable and $(4, 4)$-colorable. 

\indent Let $G$ be a graph without 4-cycles and 5-cycles. 
We show that the problem to determine whether $G$ is $(0,k)$-colorable 
is NP-complete for each positive integer $k.$  
Moreover, we construct non-$(1,k)$-colorable planar graphs 
without 4-cycles and 5-cycles for each positive integer $k.$
Finally, we prove that $G$ is $(d_1,d_2)$-colorable 
where $(d_1,d_2)=(4,4), (3,5),$ and $(2,9).$  

\section{NP-completeness of $(0,k)$-colorings} 
\begin{thm}\label{NP1}\cite{mong4}  
Let $g_{k,j}$ be the largest integer $g$ 
such that there exists a planar graph of girth $g$ that is not $(k,j)$-colorable. 
The problem to determine whether a planar graph with girth 
$g_{k,j}$ is $(k,j)$-colorable for $(k,j) \neq (0,0)$ is NP-complete. 
\end{thm}

\begin{thm}\label{NP2}
The problem to determine whether a planar graph without 4-cycles and 5-cycles 
is $(0,k)$-colorable is NP-complete for each positive integer $k.$  
\end{thm} 
\begin{pr} 
We use a reduction from the problem in Theorem \ref{NP1} 
to prove that $(0,k)$-coloring for planar graph without 4-cycles and 5-cycles.  
From \cite{mong4}, $6 \leq g_{0,1} \leq 10.$ 
Let $G$ be a graph of girth $g_{0,1}.$ 
Take $k-1$ copies of 3-cycles $v_i v'_i v''_i$ $(i=1,\ldots,k-1)$ 
for each vertex $v$ of $G.$
The graph $H_k$ is obtained from $G$ by identifying 
$v_i$ (in a 3-cycle $v_i v'_i v''_i$) to $v$ for each vertex $v.$  
The resulting graph $H_k$ has neither 4-cycles nor 5-cycles.

Suppose $G$  has a $(0,1)$-coloring $c.$ 
We extend a coloring to $c(v'_i)=1$ and $c(v''_i)=2$ 
for each vertex $v$ and each $i=1,\ldots, k-1.$ 
One can see that $c$ is a$(0,k)$-coloring of $H_k.$ 
Suppose $H_k$  has a $(0,k)$-coloring $c.$ 
Consider $v \in V(G)$ with $c(v) =2.$ 
By construction, 
$v$ has at least $k-1$ neighbors with the same color in $V(H_k) - V(G).$ 
Thus $v$ has at most one neighbor with the same color in $V(H_k) - V(G).$ 
It follows that $c$ with restriction to $V(G)$ is a $(0,1)$-coloring of $G.$ 
Hence $G$ is $(0,1)$-colorable if and only if $H_k$ is $(0,k)$-colorable. 
This completes the proof. 
\end{pr}

\section{Non-$(1,k)$-colorable planar graphs without 4-cycles and 5-cycles} 

We construct a non-$(1,k)$-colorable planar graph $G$ without 4-cycles and 5-cycles. 
Consider the graph $H_{u,v}$ shown in Figure 1. 

\begin{center}
\scalebox{0.7} 
{
\begin{pspicture}(0,-5.6259413)(19.66974,5.5859413)
\rput{-4.372639}(-0.25628725,1.0641237){\psellipse[linewidth=0.04,linestyle=dashed,dash=0.16cm 0.16cm,dimen=outer](13.808568,3.8886275)(0.112,1.5546079)}
\rput{-4.372639}(-0.2538841,1.064462){\psellipse[linewidth=0.04,linestyle=dashed,dash=0.16cm 0.16cm,dimen=outer](13.814198,3.857323)(0.216,1.5546079)}
\rput{-4.372639}(-0.25607806,1.065943){\psellipse[linewidth=0.04,linestyle=dashed,dash=0.16cm 0.16cm,dimen=outer](13.832498,3.8867977)(0.36,1.5546079)}
\rput{114.028076}(18.689152,-8.8780575){\psellipse[linewidth=0.04,linestyle=dashed,dash=0.16cm 0.16cm,dimen=outer](12.225769,1.6261551)(0.112,1.5546079)}
\rput{114.028076}(18.742256,-8.872842){\psellipse[linewidth=0.04,linestyle=dashed,dash=0.16cm 0.16cm,dimen=outer](12.2506275,1.6459968)(0.216,1.5546079)}
\rput{114.028076}(18.695421,-8.838286){\psellipse[linewidth=0.04,linestyle=dashed,dash=0.16cm 0.16cm,dimen=outer](12.215997,1.6480755)(0.36,1.5546079)}
\rput{72.36679}(12.332929,-13.307582){\psellipse[linewidth=0.04,linestyle=dashed,dash=0.16cm 0.16cm,dimen=outer](15.2632475,1.7767394)(0.112,1.5546079)}
\rput{72.36679}(12.353447,-13.339037){\psellipse[linewidth=0.04,linestyle=dashed,dash=0.16cm 0.16cm,dimen=outer](15.295009,1.7750382)(0.216,1.5546079)}
\rput{72.36679}(12.359794,-13.298567){\psellipse[linewidth=0.04,linestyle=dashed,dash=0.16cm 0.16cm,dimen=outer](15.270518,1.7996118)(0.36,1.5546079)}
\psdots[dotsize=0.2,dotangle=-0.41356933](13.930663,5.465545)
\psdots[dotsize=0.2,dotangle=-0.41356933](16.800848,1.3247203)
\psdots[dotsize=0.2,dotangle=-0.41356933](10.698409,1.00876)
\psdots[dotsize=0.2,dotangle=-0.41356933](13.707859,2.3070714)
\psline[linewidth=0.04cm](10.738843,1.0689027)(13.890954,5.506266)
\psline[linewidth=0.04cm](13.931239,5.54511)(16.781425,1.4044294)
\rput{123.139694}(12.099944,-12.8187065){\psellipse[linewidth=0.04,linestyle=dashed,dash=0.16cm 0.16cm,dimen=outer](9.519856,-3.13403)(0.112,1.5546079)}
\rput{123.139694}(12.152749,-12.800239){\psellipse[linewidth=0.04,linestyle=dashed,dash=0.16cm 0.16cm,dimen=outer](9.541261,-3.110502)(0.216,1.5546079)}
\rput{123.139694}(12.096478,-12.776638){\psellipse[linewidth=0.04,linestyle=dashed,dash=0.16cm 0.16cm,dimen=outer](9.506737,-3.1139338)(0.36,1.5546079)}
\rput{241.5404}(20.77729,6.347413){\psellipse[linewidth=0.04,linestyle=dashed,dash=0.16cm 0.16cm,dimen=outer](12.278317,-3.0118492)(0.112,1.5546079)}
\rput{241.5404}(20.724985,6.331544){\psellipse[linewidth=0.04,linestyle=dashed,dash=0.16cm 0.16cm,dimen=outer](12.247441,-3.0042126)(0.216,1.5546079)}
\rput{241.5404}(20.778954,6.306204){\psellipse[linewidth=0.04,linestyle=dashed,dash=0.16cm 0.16cm,dimen=outer](12.266881,-3.032949)(0.36,1.5546079)}
\rput{199.87912}(20.240234,2.1585915){\psellipse[linewidth=0.04,linestyle=dashed,dash=0.16cm 0.16cm,dimen=outer](10.309252,-0.69414884)(0.112,1.5546079)}
\rput{199.87912}(20.196404,2.20337){\psellipse[linewidth=0.04,linestyle=dashed,dash=0.16cm 0.16cm,dimen=outer](10.291261,-0.6679192)(0.216,1.5546079)}
\rput{199.87912}(20.199213,2.1350813){\psellipse[linewidth=0.04,linestyle=dashed,dash=0.16cm 0.16cm,dimen=outer](10.286682,-0.7023097)(0.36,1.5546079)}
\psdots[dotsize=0.2,dotangle=127.09876](8.194663,-3.9974177)
\psdots[dotsize=0.2,dotangle=127.09876](9.731509,0.8007595)
\psdots[dotsize=0.2,dotangle=127.09876](13.698108,-3.847433)
\psdots[dotsize=0.2,dotangle=127.09876](10.835708,-2.250854)
\psline[linewidth=0.04cm](13.682316,-3.8595765)(8.243078,-4.061305)
\psline[linewidth=0.04cm](8.094664,-4.017817)(9.643574,0.7644081)
\rput{-128.42137}(31.985622,9.540842){\psellipse[linewidth=0.04,linestyle=dashed,dash=0.16cm 0.16cm,dimen=outer](18.297825,-2.9571264)(0.112,1.5546079)}
\rput{-128.42137}(31.928377,9.538908){\psellipse[linewidth=0.04,linestyle=dashed,dash=0.16cm 0.16cm,dimen=outer](18.268736,-2.9442635)(0.216,1.5546079)}
\rput{-128.42137}(31.976171,9.49867){\psellipse[linewidth=0.04,linestyle=dashed,dash=0.16cm 0.16cm,dimen=outer](18.282911,-2.9759295)(0.36,1.5546079)}
\rput{-10.020658}(0.32998836,3.0061169){\psellipse[linewidth=0.04,linestyle=dashed,dash=0.16cm 0.16cm,dimen=outer](17.309433,-0.37892625)(0.112,1.5546079)}
\rput{-10.020658}(0.3355439,3.0060718){\psellipse[linewidth=0.04,linestyle=dashed,dash=0.16cm 0.16cm,dimen=outer](17.311954,-0.41063297)(0.216,1.5546079)}
\rput{-10.020658}(0.33107555,3.0101655){\psellipse[linewidth=0.04,linestyle=dashed,dash=0.16cm 0.16cm,dimen=outer](17.333067,-0.38310233)(0.36,1.5546079)}
\rput{-51.681946}(8.316359,11.211911){\psellipse[linewidth=0.04,linestyle=dashed,dash=0.16cm 0.16cm,dimen=outer](15.733523,-2.9799767)(0.112,1.5546079)}
\rput{-51.681946}(8.328966,11.187216){\psellipse[linewidth=0.04,linestyle=dashed,dash=0.16cm 0.16cm,dimen=outer](15.714332,-3.0053403)(0.216,1.5546079)}
\rput{-51.681946}(8.336786,11.216432){\psellipse[linewidth=0.04,linestyle=dashed,dash=0.16cm 0.16cm,dimen=outer](15.7484045,-2.9988067)(0.36,1.5546079)}
\psdots[dotsize=0.2,dotangle=-124.4623](19.53604,-3.9412024)
\psdots[dotsize=0.2,dotangle=-124.4623](14.498098,-4.0008893)
\psdots[dotsize=0.2,dotangle=-124.4623](17.65305,1.2322636)
\psdots[dotsize=0.2,dotangle=-124.4623](17.043795,-1.9881728)
\psline[linewidth=0.04cm](17.693836,1.1548016)(19.605604,-3.9413872)
\psline[linewidth=0.04cm](19.608046,-3.9959414)(14.581424,-4.0391383)
\psline[linewidth=0.04cm](10.865783,-2.224791)(13.722497,2.2846637)
\psline[linewidth=0.04cm](13.676769,2.2840586)(16.976768,-1.9559414)
\psdots[dotsize=0.2](0.67676944,0.24405855)
\psdots[dotsize=0.2](2.0967693,2.0640585)
\psdots[dotsize=0.2](2.0567694,-1.6159414)
\psline[linewidth=0.04cm](0.7167694,0.24405855)(2.0567694,2.0440586)
\psline[linewidth=0.04cm](2.0567694,2.0440586)(2.0567694,-1.5759414)
\psline[linewidth=0.04cm](2.0567694,-1.5759414)(0.67676944,0.26405856)
\psdots[dotsize=0.2,dotangle=-179.53052](6.396909,0.20987691)
\psdots[dotsize=0.2,dotangle=-179.53052](4.99187,-1.6216974)
\psdots[dotsize=0.2,dotangle=-179.53052](5.001715,2.058507)
\psline[linewidth=0.04cm](6.3569107,0.20954916)(5.031705,-1.6013703)
\psline[linewidth=0.04cm](5.031705,-1.6013703)(5.002043,2.0185082)
\psline[linewidth=0.04cm](5.002043,2.0185082)(6.3970733,0.18987758)
\psline[linewidth=0.04cm](2.0167694,-1.5959414)(4.9767694,-1.5959414)
\psdots[dotsize=0.2](1.9567695,-3.8959415)
\psdots[dotsize=0.2](4.9967694,-3.8959415)
\psline[linewidth=0.04cm,linestyle=dashed,dash=0.16cm 0.16cm](1.9567695,-3.8959415)(5.0367694,-3.8959415)
\usefont{T1}{ptm}{m}{n}
\rput(1.7520819,-3.4559414){\Large $u$}
\usefont{T1}{ptm}{m}{n}
\rput(5.112082,-3.4559414){\Large $v$}
\usefont{T1}{ptm}{m}{n}
\rput{-0.66679245}(-0.030569343,0.021911671){\rput(1.8421956,2.6236095){\Large $a$}}
\usefont{T1}{ptm}{m}{n}
\rput{-0.66679245}(-0.030118175,0.06124235){\rput(5.221968,2.6045077){\Large $c$}}
\usefont{T1}{ptm}{m}{n}
\rput{-0.66679245}(0.024124695,0.021244274){\rput(1.8121957,-2.0763905){\Large $b$}}
\usefont{T1}{ptm}{m}{n}
\rput{-0.66679245}(0.024807258,0.06034085){\rput(5.171968,-2.1154923){\Large $d$}}
\usefont{T1}{ptm}{m}{n}
\rput{-0.66679245}(-0.008096276,0.003974983){\rput(0.3121957,0.6836095){\Large $u$}}
\usefont{T1}{ptm}{m}{n}
\rput{-0.66679245}(-0.007677414,0.0777539){\rput(6.651968,0.6845076){\Large $v$}}
\usefont{T1}{ptm}{m}{n}
\rput(3.2813008,-5.1359415){\Huge $H_{u,v}$}
\usefont{T1}{ptm}{m}{n}
\rput(13.8113,-5.1159415){\Huge $G$}
\end{pspicture} 
}\\
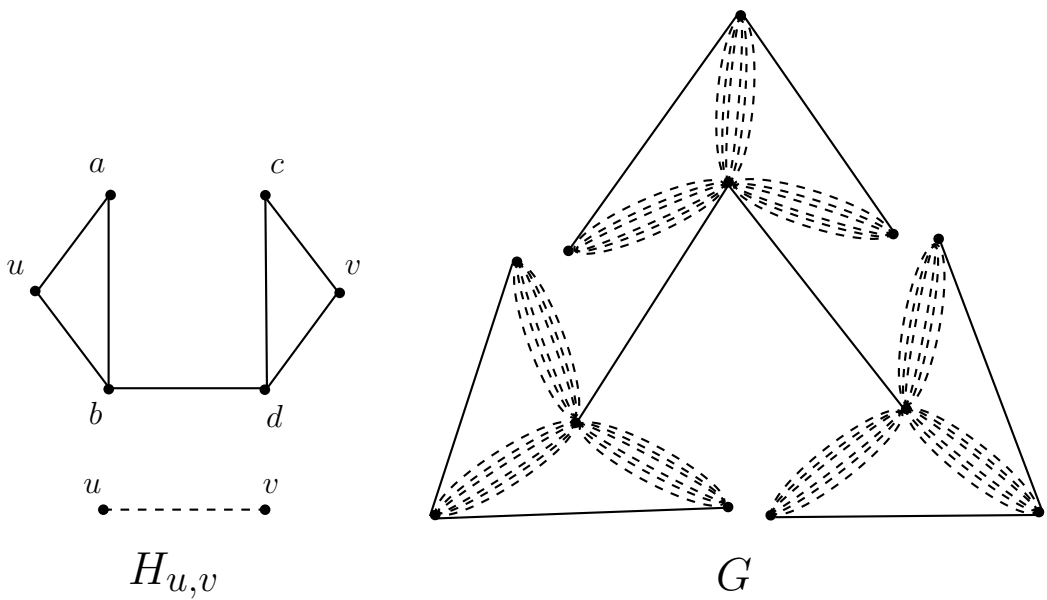
\noindent {\bf Figure 1.} 
\end{center}

A non-$(k,1)$-colorable planar graph $G$ without $4$-cycles and $5$-cycles \label{nonk1} 
The vertices $a,b,c,$ and $d$ cannot receive the same color 1. 
Now, we construct the graph $S_z$ as follows. 
Let $z$ be a vertex and $x_1 x_2 x_3$ be a path. 
Take $2k+1$ copies $H_{u_i,v_j}$ of $H_{u,v}$ with $1 \leq i \leq 2k+1$ 
and $1 \leq j \leq 3.$ 
Identify every $u_i$ with $z$ and identify $v_j$ with $x_j.$ 
Finally, we obtain $G$ from three copies $S_{z_1},S_{z_1},$ and $S_{z_3}$ 
by adding the edges $z_1z_2$ and $z_2z_3.$ 
In every $(1,k)$-coloring of $G,$ 
the path $z_1z_2z_3$ contains a vertex $z$ with color $2.$  
In the copy of $S_z$ corresponding to $z,$   
the path $x_1x_2x_3$ contains a vertex $x$ with color $2.$  
Since each of $z$ and $x$ has at most $k$ neighbors colored $2,$ 
one of $2k+1$ copies of $H_{u,v}$ between $z$ and $x,$ 
does not contain a neighbor of $z$ and $x$ colored $2.$ 
This copy is not $(1,k)$-colorable, and thus $G$ is not $(1,k)$-colorable.

\section{Helpful Tools}

\indent 
Now, we investigate $(d_1,d_2)$ such that $G$ is $(d_1,d_2)$-colorable 
for every graph $G$ without $4$-cycles and $5$-cycles.  
From two previous sections, we have that $d_1,d_2 \geq 2.$ 
First, we present useful  proposition and lemmas  
about a minimal planar graph $G$ that is not $(d_1,d_2)$-colorable where $d_1\leq d_2$.
\begin{prop}\label{l1}
(a) Each vertex $v$ of $G$  is a $2^+$-vertex.\\
(b) If $v$ is a $k$-vertex has $\alpha$ incident $3$-faces, $\beta$ adjacent good 2-vertices, and $\gamma$  pendant $3$-faces, then $\alpha\leq\lfloor \frac{k}{2} \rfloor$ and  $2\beta+\alpha+\gamma\leq k$
\end{prop}
\begin{lem}\label{l2} \cite{choi35} Let $G$ be $(d_1,d_2)$-colorable where $d_1\leq d_2$.\\
(a) If $v$ is a $3^-$-vertex, then at least two neighbors of $v$ are $(d_1+2)^+$-vertices one of which is a $(d_2+2)^+$-vertex.\\
(b) If $v$ is a $(d_1+d_2+1)^-$-vertex, then at least one neighbor of $v$ is a  $(d_1+2)^+$-vertex.
\end{lem}

\begin{lem}\label{l3} If a $2$-vertex $v$ is on a bad $3$-face $f$, then the other face $g$ which is incident to $v$ is a $7^+$-face. 
\end{lem}
\begin{pr} Suppose that a face $g$ is a $6^-$-face. Let a face $f=uvw$. By condition of $G$, a face $g$ is neither $4,5$-face nor $3$-face, otherwise $G$ contains $C_4$. Now we suppose a face $g$ is a $6$-face and let $g=u_1u_2u_3uvw$. Since $u$ is adjacent to $w$, there is a $5$-cycle $=u_1u_2u_3uw$, a contradiction. 
\end{pr}
\begin{lem}\label{l4}  Let $f$ be a $k$-face where $k\geq 7$. Then, $f$ has at most $k-6$ incident bad $2$-vertices.
\end{lem}
\begin{pr}  
By proof of Lemma \ref{l3}, if a face $f$ is incident to $m$ bad $2$-vertices, then there is a cycle $C_{k-m}$ since we can add some edge to $f$ to obtain a new cycle that has the length least than a face $f$. 
\end{pr}
\begin{lem} \label{l5}
Let $(u,v,w)$ be a bad $3$-face $f$ where $d(u)=2$.  Then at least one of  following statements is true.\\
(S1) A vertex $v$  is a $(d_1+3)^+$-vertex which has at least two $(d_2+2)$-neighbors.\\
(S2) A vertex $w$ is  a $(d_2+3)^+$-vertex which has at least two $(d_1+2)$-neighbors.\\
(S3) A vertex $v$ or a vertex $w$ is a $(d_1+d_2+2)^+$-vertex.
\end{lem}
\begin{pr} Assume $c$ is a $(d_1,d_2)$-coloring in $G-u$. If two neighbors of $u$ share the same color, then  we can color $u$ by $\{1,2\}-\{c(v)\}$. So $c(v)\not= c(w)$. By symmetry let $c(v)=1$ and $c(w)=2$.
By Lemma \ref{l2}, we have a vertex $v$  is a $(d_1+2)^+$
 and a vertex $w$  is a $(d_2+2)^+$. Then  $v$ has $d_1$ neighbors of color $1$ to forbid $u$ from being colored by $1$ and $w$ has $d_2$ neighbors of color $2$ to forbid $u$ from being colored by $2$. Next, to avoid recoloring $v$ by $2$ and $w$ by $1$. Then $v$ has one neighbor with color $2$ which has $d_2$ neighbors of color $2$ or $v$ has $d_2$ neighbors with color $2$. Otherwise,  $w$ has one neighbor with color $1$ which has $d_1$ neighbors of color $1$ or $w$ has $d_1$ neighbors with color $1$. 
\end{pr}

\section{$(4,4)$-coloring}
\begin{thm} \label{thm1} If  $G$ is a planar graph without cycles of length $4$ or $5$, then $G$ is $(4,4)$-colorable.
\end{thm}
\begin{pr}
\indent Suppose that $G$ is a minimal counterexample. The discharging process is as follows. Let the initial charge of a vertex $u$ in $G$ be $\mu(u)=2d(u)-6$ and the initial charge of a face $f$ in $G$ be $\mu(f)=d(f)-6$. Then by Euler's formula $|V(G)|-|E(G)|+|F(G)|=2$ and by the Handshaking lemma, we have
$$\displaystyle\sum_{u\in V(G)}\mu(u)+\displaystyle\sum_{f\in F(G)}\mu(f)=-12.$$
\indent Now, we establish	a new charge $\mu^*(x)$ for all $x\in V(G)\cup F(G)$ by transferring charge from one element to another and the summation of new charge $\mu^*(x)$ remains $-12$. If the final charge  $\mu^*(x)\geq 0$ for all $x\in V(G)\cup F(G)$, then we get a contradiction and the prove is completed.\\
\indent The discharging rules are\\
(R1) Every $6^+$-vertex sends charge 1 to each adjacent good $2$-vertex.\\
(R2) Every $6^+$-vertex sends charge 2 to each incident  $3$-face.\\
(R3) Every $6^+$-vertex sends charge 1 to each adjacent pendant $3$-face.\\
(R4) Every $7^+$-face sends charge 1 to each incident bad $2$-vertex.\\
(R5) Every $4$-vertex or $5$-vertex sends charge 1 to each incident $3$-face.\\
(R6) Every bad $3$-face sends charge 1 to each incident $2$-vertex.\\
\indent It remains to show that resulting $\mu^*(x)\geq 0$ for all $x\in V(G)\cup F(G)$.\\
\indent It is evident  that $\mu^*(x)=\mu(x)=0$ if $x$ is a $3$-vertex or a $6$-face.\\
\indent Now, let $v$ be a $k$-vertex.\\
\indent For $k=2$, a vertex $v$ has two $6^+$-neighbors by Lemma \ref{l2}. If $v$ is a good $2$-vertex, then $\mu^*(v)\geq \mu(v)+2\cdot1=0$ by (R1). If $v$ is a bad $2$-vertex, then $v$ is incident to a $7^+$-face by Lemma \ref{l3}. Thus $\mu^*(v)\geq \mu(v)+1+1=0$ by (R4) and (R6).\\
\indent For $k=4,5$, by Proposition \ref{l1} (b), a vertex $v$ is incident to at most two $3$-faces. By (R5),  $\mu^*(v)\geq \mu(v)-2\cdot1\geq 0$.\\
\indent Consider $k=6^+.$ Let $v$ have $\alpha$ incident $3$-faces,  $\beta$ adjacent good 2-vertices, and  $\gamma$ pendant $3$-faces. By Proposition \ref{l1} (b), we have $2\alpha+\beta+\gamma\leq d(v)$. Moreover, $\mu(v)=2d(v)-6\geq d(v)$ if $d(v)\geq 6$. Thus, by (R1), (R2), and (R3), we have $\mu^*(v)=\mu(v)-(2\alpha+\beta+\gamma)\geq  0.$\\
\indent  Now let $f$ be a $k$-face.\\
\indent For $k=7^+$, by Lemma \ref{l4}, a $k$-face $f$ has at most  $k-6$ incident   bad $2$-vertices. By (R4), $\mu^*(f)= \mu(f)-(k-6)\cdot1= 0$.\\
\indent Consider $k=3$. If $f$ is a bad $3$-face, then we have $f=(2,6^+,6^+)$-face by Lemma \ref{l2}. Then by (R2) and (R6), $\mu^*(f)\geq \mu(f)+2\cdot2-1= 0$. Now, It remains to consider a good $3$-face.
If $f$ is incident to a $4^+$-vertex and a $6^+$-vertex, then $\mu^*(f)\geq \mu(f)+2+1\geq 0$ by (R2) and (R5). If $f$ is a $(3,3,6^+)$-face, then the pendant neighbor of a $3$-vertex is a $6^+$-vertex by Lemma \ref{l2}. Thus $\mu^*(f)\geq \mu(f)+2+1+1\leq 0$ by (R2) and (R3). Finally, if $f$ is a $(4^+,4^+,4^+)$-face, then $\mu^*(f)\geq \mu(f)+3\cdot1\leq 0$ by (R5). \\
\indent  Since $\mu^*(x)\geq 0$ for all $x\in V(G)\cup F(G)$, this completes the proof.
\end{pr}

\section{$(3,5)$-coloring}
\begin{thm} If  $G$ is a planar graph without cycles of length $4$ or $5$, then $G$ is $(3,5)$-colorable.
\end{thm}
\begin{pr}
\indent Suppose that $G$ is a minimal counterexample. The discharging process is as follows. Let the initial charge of a vertex $u$ in $G$ be $\mu(u)=2d(u)-6$ and the initial charge of a face $f$ in $G$ be $\mu(f)=d(f)-6$. Then by Euler's formula $|V(G)|-|E(G)|+F(G)=2$ and by the Handshaking lemma, we have
$$\displaystyle\sum_{u\in V(G)}\mu(u)+\displaystyle\sum_{f\in F(G)}\mu(f)=-12.$$
\indent Now, we establish	a new charge $\mu^*(x)$ for all $x\in V(G)\cup F(G)$ by transferring charge from one element to another and the summation of new charge $\mu^*(x)$ remains $-12$. If the final charge  $\mu^*(x)\geq 0$ for all $x\in V(G)\cup F(G)$, then we get a contradiction and the prove is completed.\\
\indent The discharging rules are\\
(R1) Every $5$-vertex sends charge $\frac{4}{5}$ to each adjacent good $2$-vertex.\\
(R2) Every $5$-vertex sends charge $\frac{8}{5}$ to each incident  $3$-face.\\
(R3) Every $5$-vertex sends charge $\frac{4}{5}$ to each adjacent pendant $3$-face.\\
(R4) Every $6$-vertex sends charge 1 to eeach adjacent good $2$-vertex.\\
(R5) Every $6$-vertex  or $7$-vertex sends charge 2 to each incident  $3$-face.\\
(R6) Every $6$-vertex sends charge 1 to each adjacent pendant $3$-face.\\
(R7) Every $7^+$-vertex sends charge $\frac{6}{5}$ to each adjacent good $2$-vertex.\\
(R8) Every $8^+$-vertex sends charge $\frac{12}{5}$ to each incident  $3$-face.\\
(R9) Every $7^+$-vertex sends charge $\frac{6}{5}$ to each adjacent pendant $3$-face.\\
(R10) Every $7^+$-face sends charge 1 to each incident  bad $2$-vertex.\\
(R11) Every $4$-vertex  sends charge 1 to each incident $3$-face.\\
(R12) Every bad $3$-face sends charge 1 to each incident $2$-vertex.\\
\indent Next, we  show that the final charge  $\mu^*(u)$  is nonnegative.

\indent It is evident  that $\mu^*(x)=\mu(x)=0$ if $x$ is a $3$-vertex or a $6$-face.\\
\indent Now, let $v$ be a $k$-vertex.\\
\indent For $k=2$, a vertex $v$ has two $5^+$-neighbors one of which is a $7^+$-neighbor by Lemma \ref{l2}. If $v$ is a good $2$-vertex, then $\mu^*(v)\geq \mu(v)+\frac{4}{5}+\frac{6}{5}=0$ by (R1) and (R7). If $v$ is  a bad $2$-vertex, then $v$ is incident to a $7^+$-face by Lemma \ref{l3}. Thus $\mu^*(v)\geq \mu(v)+1+1=0$ by (R10) and (R12).\\
\indent For $k=4$, by Proposition \ref{l1} (b), a vertex $v$ is incident to at most two $3$-faces. By (R11),  $\mu^*(v)\geq \mu(v)-2\cdot1\geq 0$.\\
\indent Consider $k=5$. Let $v$ have $\alpha$ incident $3$-faces,  $\beta$ adjacent good 2-vertices, and  $\gamma$ pendant $3$-faces. By Proposition \ref{l1} (b), $2\alpha+\beta+\gamma\leq d(v)$.  Moreover, we have $\frac{8}{5}\alpha+\frac{4}{5}\beta+\frac{4}{5}\gamma=\frac{4}{5}(2\alpha+\beta+\gamma)\leq \frac{4}{5}d(v)$ and $\mu(v)=2d(v)-6= \frac{4}{5}d(v)$ if $d(v)=5$. Thus by (R1), (R2), and (R3), we have $\mu^*(v)=\mu(v)-(\frac{8}{5}\alpha+\frac{4}{5}\beta+\frac{4}{5}\gamma)\geq 0$.\\
\indent Consider $k=6$. Let $v$ have $\alpha$ incident $3$-faces,  $\beta$ adjacent good 2-vertices, and  $\gamma$ pendant $3$-faces. By Proposition \ref{l1} (b), we have $2\alpha+\beta+\gamma\leq d(v)$. Moreover, $\mu(v)=2d(v)-6= d(v)$ if  $d(v)= 6$. Thus, by (R4), (R5), and (R6), we have $\mu^*(v)=\mu(v)-(2\alpha+\beta+\gamma)= 0.$\\
\indent Consider $k=7$. 
If $v$ is not incident to a $3$-face, then we have $\mu^*(v)=\mu(v)-6\cdot\frac{6}{5}\geq 0$ by Lemma \ref{l2}, (R7), and (R9). If $v$ is incident to one $3$-face, then we have $\mu^*(v)=\mu(v)-(2+5\cdot\frac{6}{5})= 0$ by (R5), (R7), and (R9). If $v$ is incident to two $3$-faces, then we have $\mu^*(v)=\mu(v)-(2\cdot2+3\cdot\frac{6}{5})\geq 0$ by (R5), (R7), and (R9). Finally, if $v$ is incident to three $3$-faces, then we have $\mu^*(v)=\mu(v)-(3\cdot2+\frac{6}{5})\geq 0$ by (R5), (R7) and (R9).\\\\

\indent Consider $k=8^+$. Let $v$ have $\alpha$ incident $3$-faces,  $\beta$ adjacent good 2-vertices, and  $\gamma$ pendant $3$-faces. By Proposition \ref{l1} (b), $2\alpha+\beta+\gamma\leq d(v)$.  Moreover, we have $\frac{12}{5}\alpha+\frac{6}{5}\beta+\frac{6}{5}\gamma=\frac{6}{5}(2\alpha+\beta+\gamma)\leq \frac{6}{5}d(v)$ and $\mu(v)=2d(v)-6\geq \frac{6}{5}d(v)$ if $d(v)\geq 8$. Thus by (R7), (R8), and (R9), we have $\mu^*(v)=\mu(v)-(\frac{12}{5}\alpha+\frac{6}{5}\beta+\frac{6}{5}\gamma)\geq 0$.\\

\indent  Now let $f$ be a $k$-face.\\
\indent For, $k=7^+$. By Lemma \ref{l4}, a $k$-face $f$ has at most $k-6$ incident bad $2$-vertices. By (R11), $\mu^*(f)= \mu(f)-(k-6)\cdot1= 0$.\\
\indent Consider $k=3$. If $f$ is a bad $3$-face, then we have $f$ is a$(2,6^+,6^+)$-face 
 or $f$ is a $(2,5^+,8^+)$ by Lemma \ref{l5}. Then by (R2), (R5), (R8), and (R12), $\mu^*(f)\geq \mu(f)+2\cdot2-1= 0$ or $\mu^*(f)\geq \mu(f)+\frac{8}{5}+\frac{12}{5}-1= 0$. Now, it remains to consider a good $3$-face. If $f$ is incident to a $4^+$-vertex and a $6^+$-vertex, then $\mu^*(f)\geq \mu(f)+2+1\geq 0$ by  (R5) and (R11). If $f$ is a $(3,3,7^+)$-face, then the pendant neighbor of a $3$-vertex is a $5^+$-vertex by Lemma \ref{l2}. Thus $\mu^*(f)\geq \mu(f)+2\cdot\frac{4}{5}+2\geq 0$ by (R3) and (R5).  If $f$ is a $(3,3,5^+)$-face, then the pendant neighbor of a $3$-vertex is a $7^+$-vertex by Lemma \ref{l2}. Thus $\mu^*(f)\geq \mu(f)+2\cdot\frac{6}{5}+\frac{8}{5}\geq 0$ by (R2) and (R7).  Finally, if $f$ is a $(4^+,4^+,4^+)$-face, then $\mu^*(f)\geq \mu(f)+3\cdot1\leq 0$ by (R11).\\
\indent  Since $\mu^*(x)\geq 0$ for all $x\in V(G)\cup F(G)$, this completes the proof.
\end{pr}

\section{$(2,9)$-coloring}
\begin{thm} If  $G$ is a planar graph without cycles of length $4$ or $5$, then $G$ is $(2,9)$-colorable.
\end{thm}
\begin{pr}
\indent Suppose that $G$ is a minimal counterexample. The discharging process is as follows. Let the initial charge of a vertex $u$ in $G$ be $\mu(u)=2d(u)-6$ and the initial charge of a face $f$ in $G$ be $\mu(f)=d(f)-6$. Then by Euler's formula $|V(G)|-|E(G)|+F(G)=2$ and by the Handshaking lemma, we have
$$\displaystyle\sum_{u\in V(G)}\mu(u)+\displaystyle\sum_{f\in F(G)}\mu(f)=-12.$$
\indent Now, we establish	a new charge $\mu^*(x)$ for all $x\in V(G)\cup F(G)$ by transferring charge from one element to another and the summation of new charge $\mu^*(x)$ remains $-12$. If the final charge  $\mu^*(x)\geq 0$ for all $x\in V(G)\cup F(G)$, then we get a contradiction and the prove is completed.\\
\indent The discharging rules are\\
(R1) Every $k$-vertex for $4\leq k\leq 10$ sends charge $\frac{1}{2}$ to each adjacent good $2$-vertex.\\
(R2) Every $4$-vertex sends charge $1$ to each incident  $3$-face.\\
(R3) Every $k$-vertex for $4\leq k\leq 10$ sends $\frac{1}{2}$ to each adjacent pendant $3$-face.\\
(R4) Every $k$-vertex for $5\leq k\leq 10$ sends charge $\frac{3}{2}$ to each incident  $3$-face.\\
(R5) Every $11$-vertex sends charge $\frac{5}{2}$ to each incident  $3$-face.\\
(R6) Every $11^+$-vertex sends charge $\frac{3}{2}$ to each adjacent good $2$-vertex.\\
(R7) Every $12^+$-vertex sends charge $3$ to each incident  $3$-face.\\
(R8) Every $11^+$-vertex sends charge $\frac{3}{2}$ to each adjacent pendant $3$-face.\\
(R9) Every $7^+$-face sends charge 1 to each incident  bad $2$-vertex.\\
(R10) Every bad $3$-face sends charge 1 to each incident $2$-vertex.\\
\indent Next, we  show that the final charge  $\mu^*(u)$  is nonnegative.

\indent It is evident  that $\mu^*(x)=\mu(x)=0$ if $x$ is a $3$-vertex or a $6$-face.\\
\indent Now, let $v$ be a $k$-vertex.\\
\indent For $k=2$, a vertex $v$ has two $4^+$-neighbors one of which is a $11^+$-neighbor by Lemma \ref{l2}. If $v$ is a good $2$-vertex, then $\mu^*(v)\geq \mu(v)+\frac{1}{2}+\frac{3}{2}=0$ by (R1) and (R6). If $v$ is  a bad $2$-vertex, then $v$ is incident to a $7^+$-face by Lemma \ref{l3}. Thus $\mu^*(v)\geq \mu(v)+1+1=0$ by (R9) and (R10).\\
\indent Consider $k=4$. Let $v$ have $\alpha$ incident $3$-faces,  $\beta$ adjacent  good 2-vertices, and  $\gamma$ pendant $3$-faces. By Proposition \ref{l1} (b), $2\alpha+\beta+\gamma\leq d(v)$.  Moreover, we have $\alpha+\frac{1}{2}\beta+\frac{1}{2}\gamma=\frac{1}{2}(2\alpha+\beta+\gamma)\leq \frac{1}{2}d(v)$ and $\mu(v)=2d(v)-6= \frac{1}{2}d(v)$ if $d(v)=4$. Thus by (R1), (R2), and (R3), we have $\mu^*(v)=\mu(v)-(\alpha+\frac{1}{2}\beta+\frac{1}{2}\gamma)\geq 0$.\\
\indent Consider $k$ for $5\leq k\leq 10$. By (R1), (R3), and (R4), we show only the  case that $v$ has $\lfloor \frac{d(v)}{2} \rfloor$ incident $3$-faces because this case has final charge less than  the other cases. Consider $\frac{3}{2}\frac{d(v)}{2}\leq 2d(v)-6$, then we have $d(v)\geq 5$  because two times charge in (R1) or (R3) is less than charge in (R4). Thus we have $\mu^*(v)\geq 0.$\\
\indent Consider $k=11$. By (R5), (R6), and (R8), we  show only the case that $v$ is not incident to  $3$-face because this case has final charge less than  the other cases. 
we have $\mu^*(v)=16-10(\frac{3}{2})\geq 0$. If there is one $3$-face, then $\mu^*(v)=16-(9(\frac{3}{2})+\frac{5}{2})=0.$\\
\indent  Now let $f$ be a $k$-face.\\
\indent For $k=7^+$. By Lemma \ref{l4}, a $k$-face $f$ has at most $k-6$ incident  bad $2$-vertices. By (R9), $\mu^*(f)= \mu(f)-(k-6)\cdot1= 0$.\\
\indent Consider $k=3$. If $f$ is a bad $3$-face, then we have $f$ is a$(2,4^+,12^+)$-face 
 or $f$ is a $(2,5^+,11^+)$ by Lemma \ref{l5}. Then by (R2), (R4), (R5), and (R7), $\mu^*(f)\geq \mu(f)+1+3-1= 0$ or $\mu^*(f)\geq \mu(f)+\frac{3}{2}+\frac{5}{2}-1= 0$. Now, it remains to consider a good $3$-face. 
Consider $f$ is incident to exactly one $3$-vertex. If $f$ is not incident to a $11^+$-vertex, then pendant neighbor of a $3$-vertex is a $11^+$-vertex by Lemma \ref{l2}.
Thus $\mu^*(f)\geq \mu(f)+2\cdot\frac{1}{2}+\frac{3}{2}\geq 0$ by  (R2) and (R8).  If $f$ is incident to a $4^+$-vertex and a $11^+$-vertex, then $\mu^*(f)\geq \mu(f)+\frac{1}{2}+\frac{5}{2}\geq 0$ by  (R2) and (R5). 
If $f$ is a $(3,3,11^+)$-face, then the pendant neighbor of a $3$-vertex is a $4^+$-vertex by Lemma \ref{l2}. Thus $\mu^*(f)\geq \mu(f)+2\cdot\frac{1}{2}+\frac{5}{2}\geq 0$ by (R3) and (R5). 
 If $f$ is a $(3,3,4^+)$-face, then the pendant neighbor of a $3$-vertex is a $11^+$-vertex by Lemma \ref{l2}. Thus $\mu^*(f)\geq \mu(f)+2\cdot\frac{3}{2}+1\geq 0$ by (R2) and (R8).  
Finally, if $f$ is a $(4^+,4^+,4^+)$-face, then $\mu^*(f)\geq \mu(f)+3\cdot1\geq 0$ by (R2).\\
\indent  Since $\mu^*(x)\geq 0$ for all $x\in V(G)\cup F(G)$, this completes the proof.
\end{pr}

\end{document}